\documentclass[10pt]{article}

\usepackage{amsmath,amsthm,amsfonts,amssymb,latexsym,amscd}
\usepackage{pictex,rotating,multirow}
\begin{document}

\newtheorem*{theo}{Theorem}
\newtheorem*{pro}{Proposition}
\newtheorem*{cor}{Corollary}
\newtheorem*{lem}{Lemma}
\newtheorem{theorem}{Theorem}[section]
\newtheorem{corollary}[theorem]{Corollary}
\newtheorem{lemma}[theorem]{Lemma}
\newtheorem{proposition}[theorem]{Proposition}
\newtheorem{conjecture}[theorem]{Conjecture}
\newtheorem{definition}[theorem]{Definition}
\newtheorem{problem}[theorem]{Problem}
\newtheorem{remark}[theorem]{Remark}
\newtheorem{example}[theorem]{Example}
\newcommand{\Naturali}{{\mathbb{N}}}
\newcommand{\Reali}{{\mathbb{R}}}
\newcommand{\Complessi}{{\mathbb{C}}}
\newcommand{\Toro}{{\mathbb{T}}}
\newcommand{\Relativi}{{\mathbb{Z}}}
\newcommand{\HH}{\mathfrak H}
\newcommand{\KK}{\mathfrak K}
\newcommand{\LL}{\mathfrak L}
\newcommand{\as}{\ast_{\sigma}}
\newcommand{\tn}{\vert\hspace{-.3mm}\vert\hspace{-.3mm}\vert}
\def\mA{{\mathfrak A}}
\def\A{{\mathcal A}}
\def\mB{{\mathfrak B}}
\def\B{{\mathcal B}}
\def\C{{\mathcal C}}
\def\D{{\mathcal D}}
\def\F{{\mathcal F}}
\def\H{{\mathcal H}}
\def\J{{\mathcal J}}
\def\K{{\mathcal K}}
\def\L{{\cal L}}
\def\N{{\cal N}}
\def\M{{\cal M}}
\def\O{{\mathcal O}}
\def\P{{\cal P}}
\def\S{{\cal S}}
\def\T{{\cal T}}
\def\U{{\cal U}}
\def\W{{\cal W}}
\def\b{\lambda_B(P}
\def\j{\lambda_J(P}

\title{More Localized Automorphisms \\ of the Cuntz Algebras}
% On {\rm Aut}(\O_{n+1})

\author{Roberto Conti, Jason Kimberley and Wojciech Szyma{\'n}ski}

\date{}
\maketitle

\renewcommand{\sectionmark}[1]{}

\noindent
{\small \date{2008 August 11th}}

\vspace{7mm}
\begin{abstract}
We completely determine the localized automorphisms of the Cuntz
algebras $\O_n$ corresponding to permutation matrices in $M_n
\otimes M_n$ for $n=3$ and $n=4$. This result is obtained through a
combination of general combinatorial techniques and large scale
computer calculations. Our analysis proceeds according to the
general scheme proposed in a previous paper, where we analyzed in
detail the case of $\O_2$ using labeled rooted trees. We also
discuss those proper endomorphisms of these Cuntz algebras which
restrict to automorphisms of their respective diagonals. In the case
of $\O_3$ we compute the number of automorphisms of the diagonal
induced by permutation matrices in $M_3 \otimes M_3 \otimes M_3$.
\end{abstract}

\vfill

\vspace{3mm}

\newpage

\section{Introduction and preliminaries}

In \cite{Cun2}, Cuntz noticed that the automorphism group of $\O_n$
has a rich structure resembling that of semisimple Lie groups and
suggested an intriguing definition of the Weyl group in this
context. However, despite the fact that the Cuntz algebras $\O_n$
have been intensively studied over the last thirty years, to date
precious little was known about the structure of these Weyl groups.
In \cite{CS}, we opened a new and promising line of investigations
of this problem. We also discussed at length the case of $\O_2$
therein. In the present paper we follow it up with an analysis
focused on the cases of $\O_3$ and $\O_4$; the main result being the
complete classification of all the permutation automorphisms of
$\O_n$ for $n = 3, 4$ arising at level two (i.e. induced by a
permutation matrix in $M_n \otimes M_n$).

Until now, only few such automorphisms were known; for example,
Archbold's flip-flop automorphism of $\O_2$ \cite{Arc}
and more generally Bogolubov (permutation) automorphisms of $\O_n$.
% by reshuffling of basis
The Matsumoto-Tomiyama automorphism of $\O_4$ \cite{MaTo} was
somewhat more complicated; it was only recently recognized that it
fits into a more general pattern (see \cite{Sz} and Section
\ref{section3} below). However, all such known automorphisms were,
in some sense, isolated examples with no general or systematic
understanding of the overall situation. Finding all automorphisms
through a case-by-case examination is infeasible due to the
exceedingly large scale of the problem; so an efficient reduction
process is necessary. One could exploit the action of inner
automorphisms and Bogolubov automorphisms in this process, but this
is insufficient to significantly reduce the computation.

In \cite{CS} we discovered a powerful algorithm to construct those permutations leading to automorphisms;
surprisingly, it relies on a certain combinatorial analysis of labeled rooted trees.
This fact appears vaguely reminiscent of quantum field theory,
although our setup has nothing to do with perturbation theory.
The aforementioned reduction is a result of purely theoretical analysis of the problem and has deep theoretical implications.
However, in order to perform subsequent massive computations that emerged,
we employed the Magma \cite{Magma} computational algebra system.

In particular, as a result of these computations,
we have obtained a complete classification of automorphisms
of $\O_n$ arising from permutations of the set of multiindices $\{1,\ldots,n\}^k$
for small values of $n$ and $k$.
As a byproduct, by a similar method we can also access those endomorphisms of the
Cuntz algebra that provide automorphisms of the diagonal.

\bigskip

We now briefly describe our notation and the setup. For any integer
$n \geq 2$, the Cuntz algebra $\O_n$ is the $C^*$-algebra generated
by $n$ isometries $S_1, \ldots, S_n$ with mutually orthogonal ranges
summing up to 1. One has the unital inclusions
$$ \O_n \supset \F_n \supset \D_n \ , $$
where $\F_n$ is the UHF algebra of type $n^\infty$ and the diagonal $\D_n$ is maximal abelian in both $\F_n$ and $\O_n$.
$\F_n$ is the closure in norm of the union $\cup_{k \in {\mathbb N}} \F_n^k$ of an increasing family of matrix algebras where,
for each $k \in {\mathbb N}$, the $C^*$-subalgebra $\F_n^k$ is isomorphic to the algebra $M_{n^k}$ of $n^k \times n^k$ complex matrices.
Similarly, $\D_n$ is the norm-closure of the union of the increasing sequence of $C^*$-algebras $\D_n^k$,
each isomorphic to the diagonal matrices in $M_{n^k}$.

There is a well-known one-to-one correspondence, $u \mapsto \lambda_u$, between
$\U(\O_n)$, the group of unitary elements in $\O_n$
and ${\rm End}(\O_n)$, the semigroup of unital $*$-endomorphisms of $\O_n$,
where $\lambda_u$ is uniquely determined by $\lambda_u (S_i) = u^* S_i$, $i = 1,\ldots, n$
(here, we follow the convention in \cite{Cun2}).

Cuntz showed that the automorphisms of $\O_n$ that restrict to automorphisms of the diagonal $\D_n$ are exactly
the automorphisms induced by elements in the (unitary) normalizer
$$N_{\O_n}(\D_n) = \{z \in \U(\O_n) \ | \ z \D_n z^* = \D_n\} \ . $$
Later, Power described in detail the structure of such normalizers,
showing that any element in $N_{\O_n}(\D_n)$ is the product of a unitary in $\D_n$ and a unitary that can be written
as a finite sum of words in the $S_i$'s and their adjoints.
In particular,
$${\rm Aut}(\O_n,\D_n) = \lambda(N_{\O_n}(\D_n))^{-1}$$
and
$${\rm Aut}(\O_n,\D_n) \cap {\rm Aut}(\O_n,\F_n) = \lambda(N_{\F_n}(\D_n))^{-1}$$
where for a subset $E \in \U(\O_n)$ we define
$$\lambda(E)^{-1} = \{\lambda_u \ | u \in E\} \cap {\rm Aut}(\O_n) \ . $$

In this paper, and in \cite{CS}, we are only concerned with the structure of ${\rm Aut}(\O_n,\D_n) \cap {\rm Aut}(\O_n,\F_n)$, which gives rise (after taking quotient) to the restricted Weyl group.

Let $P_n^k$ be the group of permutations of the set $W_n^k = \{1,\ldots,n\}^k$.
Clearly, $P_n^k$ is isomorphic to ${\mathbb P}_{n^k}$, the permutation group over $n^k$ elements.
To any $\sigma \in W_n^k$ one associates a unitary $u_\sigma \in \F_n^k$ by
$$u_\sigma = \sum_{\alpha \in W_n^k} S_{\sigma(\alpha)} S_\alpha^* \ . $$
Then $\sigma \mapsto u_\sigma$ is a group isomorphism of $P_n^k$
onto its image, denoted
$\P_n^k$, that can be further identified with the group of permutation matrices in $M_{n^k}$.

Now, it follows from the above that
$$N_{\F_n}(\D_n) = N_{\O_n}(\D_n) \cap \F_n = \U(\D_n) \cdot \P_n \simeq \U(\D_n) \rtimes \P_n \ , $$
where $\P_n = \cup_k \P_n^k$, see \cite{CS}.
Thus, as Cuntz has already shown that every unitary in $\U(\D_n)$ induces an automorphism of $\O_n$,
the problem that we are facing is to determine for which permutation matrices $w \in \P_n^k$, $k=1,2,3,\ldots$,
one has $\lambda_w \in {\rm Aut}(\O_n)$.
This is exactly the point where (rooted labeled) trees come to the rescue.
For a detailed discussion, see \cite{CS}.
Throughout the next section,
we repeatedly use results from that paper.
% without further mention.

\medskip
For the reader's benefit we include the following elementary yet useful observation
valid for all $n \geq 2$.
\begin{proposition}
% Let $w \in \F_n^k$ for some $k \geq 2$.
Let $w$ be a unitary in $\O_n$.
\begin{itemize}
\item[(a)]
% If $\lambda_w \in {\rm Aut}(\O_n)$ then $\lambda_w(\F_n)=\F_n$;
If $w \in \F_n$ then
$\lambda_w \in {\rm Aut}(\O_n)$ if and only if $\lambda_w(\F_n)=\F_n$;
\item[(b)]
If $\lambda_w \in {\rm Aut}(\O_n)$ then $\lambda_w(\D_n) = \D_n$ if and only
if $w \in N_{\O_n}(\D_n)$.
\item[(c)]
If $\lambda_w(\D_n) = \D_n$ then $\lambda_w$ is an irreducible endomorphism of $\O_n$,
i.e. $\lambda_w(\O_n)' \cap \O_n = {\mathbb C}$.
\end{itemize}
\end{proposition}
\begin{proof}
(a) Necessity has been proved in \cite[Lemma 2]{Sz}. On the other
hand, $\lambda_w(\F_n) = \F_n$ implies that $\lambda_w \in {\rm
Aut}(\O_n)$ since then $w^* \in \lambda_w(\O_n)$.

(b) This is part of the statement in \cite[Proposition 1.5]{Cun2}.

(c)
Using the assumption and the fact that $\D_n$ is maximal abelian in $\O_n$ one
obtains that
\begin{align*}
\lambda_w(\O_n)' \cap \O_n \subset \lambda_w(\D_n)' \cap \O_n = \D_n ' \cap \O_n
\subset \D_n = \lambda_w(\D_n) \subset \lambda_w(\O_n)
\end{align*}
and the conclusion readily follows from $\O_n$ being simple.
\end{proof}

\medskip
As the endomorphisms of $\O_n$ (with $n \leq 4$) considered in this
paper and in \cite{CS} are all induced by unitaries $w$ in $\cup_k
\P_n^k \subset N_{\F_n}(\D_n) = N_{\O_n}(\D_n) \cap \F_n$, when they
are automorphisms they also provide, by restriction, automorphisms
of $\D_n$ and $\F_n$; when they only satisfy the weaker condition
$\lambda_w(\D_n)=\D_n$ they still act irreducibly on $\O_n$. For
example, there are four such irreducible endomorphisms of $\O_2$
corresponding to permutations in $P_2^2$, see \cite{Kaw1,CS2}.

\section{Classification results}

According to the analysis in \cite{CS}, the search for  automorphisms of $\O_n$
induced by permutations in $P_n^k$ involves the following two main steps:

\begin{itemize}
\item[(b)] finding $n$-tuples of rooted trees with vertices suitably labeled by elements of $W_n^{k-1}$,
which satisfy \cite[Lemma 4.5]{CS} (or equivalently Proposition 4.7 in {\em loc.cit.}),

\item[(d)] verifying which of the $n$-tuples satisfying (b) above also fulfill \cite[Lemma 4.10]{CS}
(or equivalently Proposition 4.11 in {\em loc.cit.}).
\end{itemize}

In turn, the solutions to condition (b) alone provide by restriction automorphisms of the diagonal $\D_n$.

\subsection{The case of $\P_3^2$}

In this case, there are only two rooted trees with three vertices.
Condition (b) can only be satisfied for the following 3-tuples of unlabeled trees:

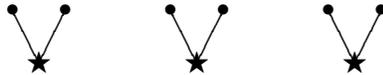
\[ \beginpicture
\setcoordinatesystem units <0.7cm,0.7cm>
\setplotarea x from 0 to 5, y from -1 to 0.1

\put {$\bullet$} at -2.5 0
\put {$\bullet$} at -3.5 0
\put {$\bigstar$} at -3 -1

\setlinear
\plot -3.5 0 -3 -1 /
\plot -2.5 0 -3 -1 /

\put {$\bullet$} at -0.5 0
\put {$\bullet$} at 0.5 0
\put {$\bigstar$} at 0 -1

\setlinear
\plot -0.5 0 0 -1 /
\plot 0.5 0 0 -1 /

\put {$\bullet$} at 2.5 0
\put {$\bullet$} at 3.5 0
\put {$\bigstar$} at 3 -1

\setlinear
\plot 2.5 0 3 -1 /
\plot 3.5 0 3 -1 /

\endpicture \]

and the three distinct 3-tuples arising by permuting

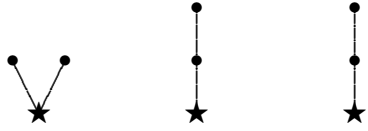
\[ \beginpicture
\setcoordinatesystem units <0.7cm,0.7cm>
\setplotarea x from 0 to 5, y from -1 to 0.1

\put {$\bullet$} at -3.5 0
\put {$\bullet$} at -2.5 0
\put {$\bigstar$} at -3 -1

\setlinear
\plot -3.5 0 -3 -1 /
\plot -2.5 0 -3 -1 /
% \put {$\alpha$} at -0.5 0
% \put {$\beta$} at -0.5 -1

\put {$\bullet$} at 0 1
\put {$\bullet$} at 0 0
\put {$\bigstar$} at 0 -1

\setlinear
\plot 0 1 0 0 /
\plot 0 0 0 -1 /

\put {$\bullet$} at 3 1
\put {$\bullet$} at 3 0
\put {$\bigstar$} at 3 -1

\setlinear
\plot 3 1 3 0 /
\plot 3 0 3 -1 /

\endpicture \]

\medskip

Thus we have four different 3-tuples. For each such 3-tuple, there
are precisely $3! \, (3!^3) = 6 \cdot 216$ permutations in $P_3^2$
satisfying condition (b) and among them $3! \, 24$ permutations
satisfying also condition (d). These figures have been obtained
through computer computations.

The corresponding labeled trees are of the form

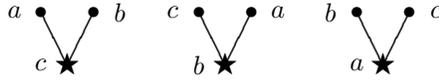
\[ \beginpicture
\setcoordinatesystem units <0.7cm,0.7cm>
\setplotarea x from 0 to 5, y from -1 to 0.1

\put {$\bullet$} at -2.5 0
\put {$\bullet$} at -3.5 0
\put {$\bigstar$} at -3 -1
\put {$c$} at -3.5 -1
\put {$a$} at -4 0
\put {$b$} at -2 0

\setlinear
\plot -3.5 0 -3 -1 /
\plot -2.5 0 -3 -1 /

\put {$\bullet$} at -0.5 0
\put {$\bullet$} at 0.5 0
\put {$\bigstar$} at 0 -1
\put {$b$} at -0.5 -1
\put {$c$} at -1 0
\put {$a$} at 1 0

\setlinear
\plot -0.5 0 0 -1 /
\plot 0.5 0 0 -1 /

\put {$\bullet$} at 2.5 0
\put {$\bullet$} at 3.5 0
\put {$\bigstar$} at 3 -1
\put {$a$} at 2.5 -1
\put {$b$} at 2 0
\put {$c$} at 4 0

\setlinear
\plot 2.5 0 3 -1 /
\plot 3.5 0 3 -1 /

\endpicture \]

and

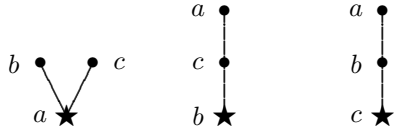
\[ \beginpicture
\setcoordinatesystem units <0.7cm,0.7cm>
\setplotarea x from 0 to 5, y from -1 to 0.1

\put {$\bullet$} at -3.5 0
\put {$\bullet$} at -2.5 0
\put {$\bigstar$} at -3 -1
\put {$a$} at -3.5 -1
\put {$b$} at -4 0
\put {$c$} at -2 0

\setlinear
\plot -3.5 0 -3 -1 /
\plot -2.5 0 -3 -1 /

\put {$\bullet$} at 0 1
\put {$\bullet$} at 0 0
\put {$\bigstar$} at 0 -1
\put {$b$} at -0.5 -1
\put {$c$} at -0.5 0
\put {$a$} at -0.5 1

\setlinear
\plot 0 1 0 0 /
\plot 0 0 0 -1 /

\put {$\bullet$} at 3 1
\put {$\bullet$} at 3 0
\put {$\bigstar$} at 3 -1
\put {$c$} at 2.5 -1
\put {$b$} at 2.5 0
\put {$a$} at 2.5 1

\setlinear
\plot 3 1 3 0 /
\plot 3 0 3 -1 /

\endpicture \]
where $a$, $b$ and $c$ are distinct elements in $\{1,2,3\}$.

In particular, for each fixed set of labels on a 3-tuple (and there are $3!$ of them) there are $24$ permutations
satisfying both the conditions (b) and (d).

\begin{example}\rm
Bogolubov automorphisms always give rise to 3-tuples of the first type \cite[Example 4.4]{CS}.
An example of a 3-tuple of the second type (with labels $a=1, b=2, c=3$)
is provided by the transposition (2,3), where we identify elements of
$W_2^2 = \{1,2,3\}^2 = \{11, 21, 31, 12, 22, 32, 13, 23, 33\}$ with \{1,2,\ldots,9\}.
\end{example}

All in all, we see that the numbers of automorphisms arising from $P_2^3$ are as follows.

\begin{theorem}
One has
\begin{align*}
 & \# \{\sigma \in P_3^2 \ : \ \lambda_{u_\sigma}|_{\D_3} \in {\rm Aut}(\D_3) \}
 = 4 \cdot 3! \cdot 216 = 5\:184, \\
 & \# \{\sigma \in P_3^2 \ : \ \lambda_{u_\sigma} \in {\rm Aut}(\O_3) \} =
4 \cdot 3! \cdot 24 = 576.
\end{align*}
In particular, there are $4\cdot24=96$ distinct classes of automorphisms in ${\rm Out}({\mathcal O}_3)$
corresponding to permutations in $P_3^2$.
\end{theorem}

The latter number has been independently verified by solving the equations in \cite[Subsection 6.1]{CS}.

\subsection{The case of $\P_3^3$}

In this case, there are 286 rooted trees with $n^{k-1}=9$ vertices,
of which 171 satisfy our basic conditions:
that each vertex has in-degree at most $n=3$
(recall that there is a loop at the root, adding 1 to its in-degree).
Let us define the \emph{in-degree type} of a rooted tree to be the multiset of the in-degrees of its vertices.

We list the 171 rooted trees in Table \ref{classified171};
they are classified by the eleven in-degree types $\{A\ldots K\}$ listed in Table \ref{in_deg_types_table}.

\begin{table}[hb]
\caption{The in-degree types for $\P_3^3$.}
\label{in_deg_types_table}
\[
\begin{tabular}{|c|rrrr|r|}
\hline
\multirow{3}{*}
{
\begin{sideways}
Type
\end{sideways}
}
&\multicolumn{4}{c|}{\footnotesize In-Degree} &
\multirow{3}{*}{
\begin{sideways}
\# trees
\end{sideways}
}\\
%\cline{2-5}
 & 0 & 1 & 2 & 3 & \\
\cline{2-5}
%$j$
& \multicolumn{4}{c|}{\scriptsize Multiplicities:}
&
\\
\hline
\hline
$A$ &
6 & 0 & 0 & 3
& 2
\\
$B$ &
5 & 1 & 1 & 2
& 18
\\
$C$ &
5 & 0 & 3 & 1
& 8
\\
$D$ &
4 & 3 & 0 & 2
& 14
\\
$E$ &
4 & 2 & 2 & 1
& 46
\\
$F$ &
4 & 1 & 4 & 0
& 9
\\
$G$ &
3 & 4 & 1 & 1
& 33
\\
$H$ &
3 & 3 & 3 & 0
& 24
\\
$I$ &
2 & 6 & 0 & 1
& 4
\\
$J$ &
2 & 5 & 2 & 0
& 12
\\
$K$ &
1 & 7 & 1 & 0
& 1
\\
\hline
\end{tabular}
\]
\end{table}

\def\jskScale{1.25mm}
\def\jskRoot{\cdot}

\begin{table}
\caption{In-degree types and trees.}
\bigskip
\label{classified171}
\begin{tabular}{|r|cccccccccccc|} %|rrrr
\hline
&\multicolumn{12}{c}{Trees of given in-degree type }\vline
\\
\hline
\hline

$A$ &
%{\small 6} & {\small 0} & {\small 0} & {\small 3} &
\beginpicture
\setcoordinatesystem units <\jskScale,\jskScale>
\put {$\cdot$} at 0 1
\plot 0 1 0 0 /
\put {$\cdot$} at 0 2
\plot 0 2 0 1 /
\put {$\cdot$} at 0 3
\plot 0 3 0 2 /
\put {$\cdot$} at 1 3
\plot 1 3 0 2 /
\put {$\cdot$} at 2 3
\plot 2 3 0 2 /
\put {$\cdot$} at 3 2
\plot 3 2 0 1 /
\put {$\cdot$} at 4 2
\plot 4 2 0 1 /
\put {$\cdot$} at 5 1
\plot 5 1 0 0 /
\put {$\jskRoot$} at 0 0
\endpicture
&
\beginpicture
\setcoordinatesystem units <\jskScale,\jskScale>
\put {$\cdot$} at 0 1
\plot 0 1 0 0 /
\put {$\cdot$} at 0 2
\plot 0 2 0 1 /
\put {$\cdot$} at 1 2
\plot 1 2 0 1 /
\put {$\cdot$} at 2 2
\plot 2 2 0 1 /
\put {$\cdot$} at 3 1
\plot 3 1 0 0 /
\put {$\cdot$} at 3 2
\plot 3 2 3 1 /
\put {$\cdot$} at 4 2
\plot 4 2 3 1 /
\put {$\cdot$} at 5 2
\plot 5 2 3 1 /
\put {$\jskRoot$} at 0 0
\endpicture
&&&&&&&&&&
\\
\hline
\end{tabular}
\end{table}

We wish to find 3-tuples $f$ of labeled trees such that
\[
    \forall{j\in\left\{1..9\right\}}\quad{f_1(j)+f_2(j)+f_3(j)=3}.
\]
For such an $f$ we define the in-degree \emph{alignment} matrix $M$ where
$M_{ij}$ is the in-degree of the vertex labeled $j$ in the tree $f_i$.
Every in-degree alignment matrix has each row adding to $n^{k-1}=9$ and each column adding to $n=3$.
In order to find all required $f$ we first determine the possible in-degree alignments of our 11 types.

Now the number of size three multisets with elements chosen from a
set of eleven is 286 (the eleventh tetrahedral number)\footnote {Is
it just a combinatorial accident that the eleventh tetrahedral
number is the same as the number of all (unordered) rooted trees on
nine-vertices.}. Of the 286 size three multisets of in-degree types,
we compute that 100 have at least one alignment. The number of
alignments (up to consistent relabeling) is 133.

After about 200 processor days we report that condition (b) is satisfied
for a set $\F$ of $7\:390$ three-tuples of labeled trees,
up to permutation of tree position (action of $S_3$) and consistent relabeling of all trees (action of $S_9$).
Only 110 of the 171 unlabeled trees appear in $\F$
(those which do not appear in $\F$ are marked in Table \ref{classified171} with a dotted backslash);
they have the first eight, $\{A\ldots H\}$, of the eleven in-degree types.
In these $\F$ there occur 474 multisets of three unlabeled trees;
they have the six distinct three element multisets of in-degree types listed in the first column of Table \ref{six_3ems_idt}.
The second column contains the number of three element multisets of unlabeled trees having the respective types;
the third column, $\F_{\mathrm{types}}$, is the partition of $\F$ according to the respective types;
the fourth column of each row in Table \ref{six_3ems_idt},
$\#f$ covered,
is the inner product of the last two columns of the corresponding table in the Appendix.

\begin{table}[h]
\caption{Three element multisets of in-degree types.}
\label{six_3ems_idt}
\[
\begin{tabular}{|ccc|r|r|r@{$\;\cdot\:9!\ $}|}
\hline
\multicolumn{3}{|c|}{\small ID types} &
\small{\# tree triples} &
%{\footnotesize$\#\left\{\left\{*t_1,t_2,t_3*\right\}\right\}$}&
\# $\F_{\mathrm{types}}$ &
\multicolumn{1}{r|}{\#$f$ \small{covered}}
\\
%&&&&& \multicolumn{1}{r|}{\footnotesize{by types}}
%\\
\hline
\hline
$A$ & $A$ & $A$ & 4 & 2 168 & 12924 \\
$A$ & $B$ & $B$ & 176 & 2 782 & 16650 \\
$A$ & $C$ & $D$ & 75 & 950 & 5700 \\
$A$ & $E$ & $E$ & 180 & 1 072 & 6396 \\
$A$ & $F$ & $G$ & 31 & 392 & 2352 \\
$A$ & $H$ & $H$ & 8 & 26 & 150 \\
%
%$A$ & $A$ & $A$ & 2\:168 & 12\:924 \\
%$A$ & $B$ & $B$ & 2\:782 & 16\:650 \\
%$A$ & $C$ & $D$ & 950 & 5\:700 \\
%$A$ & $E$ & $E$ & 1\:072 & 6\:396 \\
%$A$ & $F$ & $G$ & 392 & 2\:352 \\
%$A$ & $H$ & $H$ & 26 & 150 \\
\hline
\multicolumn{3}{|c|}{Total} & 474 & 7\:390 & 44\:172 \\
\hline
\end{tabular}
\]
\end{table}

In total, we have
\[
%           \left(1\times1+ 2\times1 + 3\times53+ 6\times7\:335
%           \right)
%       \times
%           9!
%   =
            44\:172
        \times
            9!
    =
        16\:029\:135\:360
\]
three-tuples of labeled trees satisfying condition (b).
Therefore we have
\[
%   16\:029\:135\:360 \times 6^9
            44\:172
        \times
            9!
        \times
            n!^{n^{k-1}}
    =
            44\:172
        \times
            9!
        \times
            6^9
    =
        161\:536\:753\:300\:930\:560
\]
permutations $\sigma\in P_3^3$ satisfying condition (b).

Unfortunately, at this stage we cannot provide the precise number of permutations
satisfying condition (d) as this job exceeds our computational resources:
it would take about 32 processor years to compute.

Some examples of 3-tuples of labeled trees satisfying condition (b)
are listed in the first column of Table \ref{tableP33};
the second column contains the size of the orbit of the combined actions of $S_3$ and $S_9$ on the first entry of each row;
the third column is a count of the number of permutations corresponding to $f$ that satisfy condition (d);
the last column contains (when one exists) an example permutation satisfying condition (d)
with labels $(a,b,c,\ldots,i)$ chosen to be $(1,2,3,\ldots,9)=(\,(1,1),\,(2,1),\,(3,1),\,\ldots,(3,3)\,)$.

\newcommand{\tAAABBBCCC}[9]
{
\beginpicture
\setcoordinatesystem units <2.5mm,2.5mm>
%\setplotarea x from 0 to 2, y from 0 to 3
\put {\tiny$#1$} at 6.5 -0.5
\put {\tiny$#2$} at 4.5 1
\put {\tiny$#3$} at 8.5 1
\put {\tiny$#4$} at 4  2.5
\put {\tiny$#5$} at 5  2.5
\put {\tiny$#6$} at 6  2.5
\put {\tiny$#7$} at 7  2.5
\put {\tiny$#8$} at 8  2.5
\put {\tiny$#9$} at 9  2.5
\put {$\cdot$} at 6.5 0
\put {$\cdot$} at 5   1
\put {$\cdot$} at 8   1
\put {$\cdot$} at 4   2
\put {$\cdot$} at 5   2
\put {$\cdot$} at 6   2
\put {$\cdot$} at 7   2
\put {$\cdot$} at 8   2
\put {$\cdot$} at 9   2
\setlinear
\plot 5 1 6.5 0 /
\plot 8 1 6.5 0 /
\plot 4 2 5 1 /
\plot 5 2 5 1 /
\plot 6 2 5 1 /
\plot 7 2 8 1 /
\plot 8 2 8 1 /
\plot 9 2 8 1 /
\endpicture
}

\newcommand{\tAABCCCBBA}[9]
{
\beginpicture
\setcoordinatesystem units <2.5mm,2.5mm>
%\setplotarea x from 0 to 2, y from 0 to 3
\put {\tiny$#1$} at  2.6 -0.4
\put {\tiny$#2$} at  1.6  0.6
\put {\tiny$#3$} at  0.6  1.6
\put {\tiny$#4$} at -0.4  2.6
\put {\tiny$#5$} at 1   3.5
\put {\tiny$#6$} at 2.4 3.4
\put {\tiny$#7$} at 2.4 2.4
\put {\tiny$#8$} at 3.4 2.4
\put {\tiny$#9$} at 4.4 1.4
\put {$\cdot$} at 3 0
\put {$\cdot$} at 2 1
\put {$\cdot$} at 1 2
\put {$\cdot$} at 0 3
\put {$\cdot$} at 1 3
\put {$\cdot$} at 2 3
\put {$\cdot$} at 2 2
\put {$\cdot$} at 3 2
\put {$\cdot$} at 4 1
\setlinear
\plot 0 3 3 0 /
\plot 1 3 1 2 /
\plot 2 3 1 2 /
\plot 2 2 2 1 /
\plot 3 2 2 1 /
\plot 4 1 3 0 /
\endpicture
}

\newcommand{\tECCBDDEFF}[9]
{
\beginpicture
\setcoordinatesystem units <2.5mm,2.5mm>
\setplotarea x from 1 to 6, y from 0 to 4
\put {\tiny$#3$} at 3   0
\put {\tiny$#2$} at 3   1
\put {\tiny$#4$} at 3.5 2.7
\put {\tiny$#5$} at 1.7 2.7
\put {\tiny$#1$} at 0.6 4.3
\put {\tiny$#7$} at 2.6 4.3
\put {\tiny$#8$} at 4.4 4.3
\put {\tiny$#6$} at 5.3   2.7
\put {\tiny$#9$} at 6.4 4.3
\put {$\cdot$} at 3.5 0
\put {$\cdot$} at 3.5 1
\put {$\cdot$} at 3.5 2
\put {$\cdot$} at 2   3
\put {$\cdot$} at 1   4
\put {$\cdot$} at 3   4
\put {$\cdot$} at 4   4
\put {$\cdot$} at 5   3
\put {$\cdot$} at 6   4
\setlinear
\plot 1 4 2 3 3.5 2 3.5 0 /
\plot 3 4 2 3 /
\plot 4 4 5 3 3.5 2 /
\plot 6 4 5 3 /
\endpicture
}

\newcommand{\tAABCDEFFF}[9]
{
\beginpicture
\setcoordinatesystem units <2.5mm,2.5mm>
%\setplotarea x from 1 to 6, y from 0 to 4
\put {\tiny$#1$} at  0.5 0
\put {\tiny$#2$} at  0.5 1
\put {\tiny$#3$} at  0.5 2
\put {\tiny$#4$} at  0.5 3
\put {\tiny$#5$} at  0.5 4
\put {\tiny$#6$} at  0.4 4.8
\put {\tiny$#7$} at -0.5 6
\put {\tiny$#8$} at  1   6.5
\put {\tiny$#9$} at  2.5 6
\put {$\cdot$} at 1 0
\put {$\cdot$} at 1 1
\put {$\cdot$} at 1 2
\put {$\cdot$} at 1 3
\put {$\cdot$} at 1 4
\put {$\cdot$} at 1 5
\put {$\cdot$} at 0 6
\put {$\cdot$} at 1 6
\put {$\cdot$} at 2 6
\setlinear
\plot 1 6 1 0 /
\plot 0 6 1 5 2 6 /
\endpicture
}

\newcommand{\tAABCDEFEE}[9]
{
\beginpicture
\setcoordinatesystem units <2.5mm,2.5mm>
%\setplotarea x from -1 to 2, y from 0 to 6
\put {\tiny$#1$} at  0.5 0
\put {\tiny$#2$} at  0.5 1
\put {\tiny$#3$} at  0.5 2
\put {\tiny$#4$} at  0.5 3
\put {\tiny$#5$} at  0.5 3.9
\put {\tiny$#6$} at -0.6 4.6
\put {\tiny$#7$} at -1.5 6
\put {\tiny$#8$} at  1   5.6
\put {\tiny$#9$} at  2.5 5
\put {$\cdot$} at 1 0
\put {$\cdot$} at 1 1
\put {$\cdot$} at 1 2
\put {$\cdot$} at 1 3
\put {$\cdot$} at 1 4
\put {$\cdot$} at 0 5
\put {$\cdot$} at -1 6
\put {$\cdot$} at 1 5
\put {$\cdot$} at 2 5
\setlinear
\plot  1 5 1 0 /
\plot -1 6 0 5 1 4 2 5 /
\endpicture
}

\begin{table}
\caption{Examples for $\P_3^3$}
\label{tableP33}
\[
\begin{tabular}{|ccc|r|r|c|}
\hline
\multicolumn{3}{|c|}{Labeled Trees} &
{\# (b)} &
\# (d) &
{Example} \\
\hline
\hline
%\multicolumn{1}{c}{}\\
\hline
$A$ & $A$ & $A$ & \multicolumn{3}{c|}{} \\
\hline
\hline
%2 & 2 & 2 & \multicolumn{3}{c|}{} \\
%\hline
%  \begin{tabular}{c}
       \tAAABBBCCC{g}{h}{i}{a}{b}{c}{d}{e}{f}  &
%       8, 8, 8, 9, 9, 9, 7, 7, 7  \\
       \tAAABBBCCC{d}{e}{f}{a}{b}{c}{g}{h}{i}  &
%       5, 5, 5, 4, 4, 4, 6, 6, 6  \\
       \tAAABBBCCC{a}{b}{c}{d}{e}{f}{g}{h}{i}  &
%       1, 1, 1, 2, 2, 2, 3, 3, 3  \\
%  \end{tabular} &&&
$1\cdot 9!$ & 312 &
\begin{small}
\begin{tabular}{c}
(1, 6, 26, 7, 22, 17)\\
(2, 12, 24, 20, 18, 13, 14)\\
(3, 27, 16, 25, 19, 9, 10)\\
(4, 11, 15, 23, 8)\\
\end{tabular}
\end{small}
 \\ \hline
\hline

%  \begin{tabular}{c}
        \tAABCCCBBA{g}{h}{i}{a}{b}{c}{d}{e}{f}   &
%       9, 9, 9, 8, 8, 7, 7, 7, 8  \\
       \tAAABBBCCC{d}{e}{f}{a}{b}{c}{g}{h}{i}  &
%       5, 5, 5, 4, 4, 4, 6, 6, 6  \\
       \tAAABBBCCC{a}{b}{c}{d}{e}{f}{g}{h}{i}  &
%       1, 1, 1, 2, 2, 2, 3, 3, 3  \\
%  \end{tabular} &&&
$6\cdot 9!$ &  0 & \\ \hline

%  \begin{tabular}{c}
     \tAABCCCBBA{a}{b}{c}{d}{e}{f}{g}{h}{i}   &
%         1, 1, 2, 3, 3, 3, 2, 2, 1 \\
       \tAAABBBCCC{d}{e}{f}{a}{b}{c}{g}{h}{i}  &
%         5, 5, 5, 4, 4, 4, 6, 6, 6 \\
       \tAAABBBCCC{g}{h}{i}{a}{b}{c}{d}{e}{f}  &
%         8, 8, 8, 9, 9, 9, 7, 7, 7 \\
%  \end{tabular} &&&
$6\cdot 9!$ & 240 &
\begin{small}
\begin{tabular}{c}
(1, 25, 24, 23, 2, 19)\\
(3, 16, 27, 15, 26)\\
(4, 17, 9, 18, 12, 10)\\
(6, 20, 22, 14, 8)\\
(7, 21, 13, 11)
\end{tabular}
\end{small}
\\ \hline
\hline

%  \begin{tabular}{c}
        \tAABCCCBBA{g}{h}{i}{a}{b}{c}{d}{e}{f}   &
%       9, 9, 9, 8, 8, 7, 7, 7, 8  \\
        \tAABCCCBBA{d}{e}{f}{a}{b}{c}{g}{h}{i}   &
%       6, 6, 6, 4, 4, 5, 5, 5, 4  \\
       \tAAABBBCCC{a}{b}{c}{d}{e}{f}{g}{h}{i}  &
%       1, 1, 1, 2, 2, 2, 3, 3, 3  \\
%  \end{tabular} &&&
$3\cdot 9!$ & 216 &
\begin{small}
\begin{tabular}{c}
(1, 3, 27, 4, 26, 10, 9)\\
(2, 18, 7, 16, 19, 6)\\
(5, 20, 12, 21, 24)\\
(8, 25, 22, 11, 15)\\
(13, 14, 17)\\
\end{tabular}
\end{small}
 \\ \hline

%  \begin{tabular}{c}
        \tAABCCCBBA{g}{h}{i}{a}{b}{c}{d}{e}{f}   &
%         9, 9, 9, 8, 8, 7, 7, 7, 8 \\
        \tAABCCCBBA{a}{b}{c}{d}{e}{f}{g}{h}{i}   &
%         1, 1, 2, 3, 3, 3, 2, 2, 1 \\
       \tAAABBBCCC{d}{e}{f}{a}{b}{c}{g}{h}{i}  &
%         5, 5, 5, 4, 4, 4, 6, 6, 6 \\
%  \end{tabular} &&&
$6\cdot 9!$ & 0 & \\ \hline

%  \begin{tabular}{c}
        \tAABCCCBBA{d}{e}{f}{a}{b}{c}{g}{h}{i}   &
%         6, 6, 6, 4, 4, 5, 5, 5, 4 \\
        \tAABCCCBBA{a}{b}{c}{d}{e}{f}{g}{h}{i}   &
%         1, 1, 2, 3, 3, 3, 2, 2, 1 \\
       \tAAABBBCCC{g}{h}{i}{a}{b}{c}{d}{e}{f}  &
%         8, 8, 8, 9, 9, 9, 7, 7, 7 \\
%  \end{tabular} &&&
$6\cdot 9!$ & 0 & \\ \hline
\hline

%  \begin{tabular}{c}
        \tAABCCCBBA{g}{h}{i}{a}{b}{c}{d}{e}{f}   &
%     9, 9, 9, 8, 8, 7, 7, 7, 8  \\
        \tAABCCCBBA{d}{e}{f}{a}{b}{c}{g}{h}{i}   &
%     6, 6, 6, 4, 4, 5, 5, 5, 4  \\
        \tAABCCCBBA{a}{b}{c}{d}{e}{f}{g}{h}{i}   &
%     1, 1, 2, 3, 3, 3, 2, 2, 1  \\
%  \end{tabular} &&&
$6\cdot 9!$ &  0 & \\ \hline
\hline
\hline
$A$ & $F$ & $G$ & \multicolumn{3}{c|}{} \\
\hline
\hline
%  \begin{tabular}{c}
       \tAAABBBCCC{g}{h}{i}{a}{b}{c}{d}{e}{f}  &
%       8, 8, 8, 9, 9, 9, 7, 7, 7  \\
         \tECCBDDEFF{a}{b}{c}{d}{e}{f}{g}{h}{i}  &
%       5, 3, 3, 2, 4, 4, 5, 6, 6  \\
         \tAABCDEFFF{b}{c}{d}{e}{f}{a}{g}{h}{i}  &
%       6, 2, 2, 3, 4, 5, 1, 1, 1  \\
%  \end{tabular} &&&
$6\cdot 9!$ &  0  & \\ \hline
\hline

%  \begin{tabular}{c}
       \tAAABBBCCC{g}{h}{i}{a}{b}{c}{d}{e}{f}  &
%       8, 8, 8, 9, 9, 9, 7, 7, 7  \\
         \tECCBDDEFF{a}{b}{c}{d}{e}{f}{g}{h}{i}  &
%       5, 3, 3, 2, 4, 4, 5, 6, 6  \\
         \tAABCDEFEE{b}{c}{d}{e}{a}{f}{g}{h}{i}  &
%       5, 2, 2, 3, 4, 1, 6, 1, 1  \\
%  \end{tabular} &&&
$6\cdot 9!$ &  0 & \\ \hline
\hline
\end{tabular}
\]
\end{table}

\clearpage
\subsection{The case of $\P_4^2$}
In this case, there are four (unlabeled) rooted trees with four vertices, namely

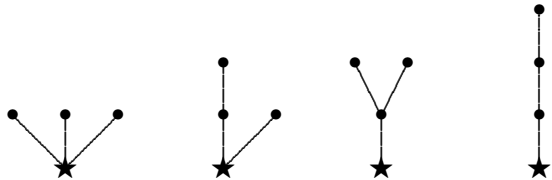
\[ \beginpicture
\setcoordinatesystem units <0.7cm,0.7cm>
\setplotarea x from 0 to 5, y from -1 to 0.1

\put {$\bullet$} at -4 0
\put {$\bullet$} at -3 0
\put {$\bullet$} at -2 0
\put {$\bigstar$} at -3 -1

\setlinear
\plot -4 0 -3 -1 /
\plot -3 0 -3 -1 /
\plot -2 0 -3 -1 /

\put {$\bullet$} at 0 1
\put {$\bullet$} at 0 0
\put {$\bullet$} at 1 0
\put {$\bigstar$} at 0 -1

\setlinear
\plot 0 1 0 0 /
\plot 1 0 0 -1 /
\plot 0 0 0 -1 /

\put {$\bullet$} at 3.5 1
\put {$\bullet$} at 2.5 1
\put {$\bullet$} at 3 0
\put {$\bigstar$} at 3 -1

\setlinear
\plot 3.5 1 3 0 /
\plot 2.5 1 3 0 /
\plot 3 0 3 -1 /
% \put {$\beta$} at 4.5 0
% \put {$\alpha$} at 4.5 -1

\put {$\bullet$} at 6 2
\put {$\bullet$} at 6 1
\put {$\bullet$} at 6 0
\put {$\bigstar$} at 6 -1

\setlinear
\plot 6 2 6 1 /
\plot 6 1 6 0 /
\plot 6 0 6 -1 /

\endpicture \]

\medskip

One verifies that only eight types of 4-tuples of such trees admit labelings satisfying condition (b).
By a type we mean an unordered set of four trees making up a 4-tuple
(two different 4-tuples belong to the same type if one can be obtained from the other by a permutation of the unlabeled trees).
These types are listed in the first column of Table \ref{table1P42} below.
The second column of this table gives the number of distinct labelings satisfying condition (b) and corresponding to each type.
These numbers are factorized as $ X \cdot Y \cdot Z$, where $X$ is the number of distinct 4-tuples of unlabeled trees of the given type,
$Y = 4!$ is the number of permutations of labels (it corresponds to action of inner automorphisms arising from $P_4^1$ \cite[Subsection 4.2]{CS}),
and $Z$ is the number of orbits under this action.
The last column contains the number of all permutations in $P_4^2$ satisfying both conditions (b) and (d) whose corresponding trees
are of the given type.

The number of permutations satisfying condition (d) depends both on
the type of the corresponding 4-tuple of trees and on the specific
labeling. However, as it turns out, it does not depend on the
permutation of unlabeled trees within the type. Precise information
to this effect is provided in Table \ref{table2P42} below. There is
a natural action of $S_4 \times S_4$ on 4-tuples of labeled trees
with four vertices, by permutation of the labels (simultaneously on
all four trees) and permutation of the four trees. Labelings
satisfying condition (b) give rise to 19 orbits for this action, and
representatives of these 19 orbits are listed in column 1. They are
further grouped according to their types. The second column
describes the partition of each orbit of the $S_4 \times S_4$-action
into orbits of an action of $S_4$ by permutation of labels. For
example, $144 = 6 \times 4!$ indicates that the corresponding $S_4
\times S_4$-orbit has 144 elements, partitioned into 6 $S_4$-orbits
with 4! = 24 elements each. The total number of permutations
satisfying condition (b) corresponding to each row is thus obtained
by multiplying the number in the second column by the combinatorial
factor $4!^{4} = 331\:776$ \cite[Section 4.2]{CS}. The third column
contains the number of permutations satisfying condition (d) for
each element in the $S_4 \times S_4$-orbit. The total number of
permutations satisfying condition (d) corresponding to a given row
is thus obtained by multiplying the numbers in the second and in the
third column. The last column of Table \ref{table2P42} contains an
example of a permutation satisfying condition (d) (if it exists)
with the choice of labels $a=1$, $b=2$, $c=3$, $d=4$ and order of
trees as given in the table.

As a consequence of the above, we obtain the following result.

\begin{theorem} One has
\begin{align*}
 & \# \{\sigma \in P_4^2 \ : \ \lambda_{u_\sigma}|_{\D_4} \in {\rm Aut}(\D_4) \} = 5\:400 \cdot 4!^4 = 1\:791\:590\:400, \\
 & \# \{\sigma \in P_4^2 \ : \ \lambda_{u_\sigma} \in {\rm Aut}(\O_4) \} = 5\:771\:520.
\end{align*}
In particular, there are $240\:480$ distinct classes of automorphisms in ${\rm Out}({\mathcal O}_4)$
corresponding to permutations in $P_4^2$.
\end{theorem}

\begin{table}
\caption{Tree types for $P_4^2$}
\label{table1P42}
\[
\begin{tabular}{|l||l|r|}
\hline
type & $\#$ (b) & $\# \sigma$ (d)\\
\hline\hline
$\alpha\alpha\alpha\alpha$ & $24 = 1 \cdot 24 \cdot 1$  & 51\:840 \\
\hline
$\alpha\alpha\beta\beta$ &  $576 = 6 \cdot 24 \cdot 4$ & 787\:968 \\
\hline
$\alpha\alpha\gamma\gamma$ & $288 = 6 \cdot 24 \cdot 2$ & 311\:040 \\
\hline
$\alpha\beta\beta\beta$ & $768 = 4 \cdot 24 \cdot 8$ & 746\:496  \\
\hline
$\alpha\beta\beta\delta$ & $1\:152 = 12 \cdot 24 \cdot 4$ & 1\:575\:936 \\
\hline
$\alpha\beta\gamma\delta$ & $1\:152 = 24 \cdot 24 \cdot 2$ & 1\:244\:160  \\
\hline
$\beta\beta\gamma\gamma$ & $1\:152 = 6 \cdot 24 \cdot 8$ & 787\:968  \\
\hline
$\gamma\gamma\gamma\gamma$ & $288 = 1 \cdot 24 \cdot 12$ & 266\:112  \\
\hline
\hline
total &5\:400& 5\:771\:520 \\
\hline
\end{tabular}
\]
\end{table}

\newcommand{\tW}[4]
{
\beginpicture
\setcoordinatesystem units <2.5mm,2.5mm>
%\setplotarea x from 0 to 2, y from 0 to 3
\put {\tiny$#1$} at 1   0.1
\put {\tiny$#2$} at 0   2.1
\put {\tiny$#3$} at 1   2.1
\put {\tiny$#4$} at 2   2.1
\put {$\cdot$} at   1   0.5
\put {$\cdot$} at   0   1.5
\put {$\cdot$} at   1   1.5
\put {$\cdot$} at   2   1.5
\setlinear
\plot 1 0.5 0 1.5 /
\plot 1 0.5 1 1.5 /
\plot 1 0.5 2 1.5 /
\endpicture
}

\newcommand{\tL}[4]
{
\beginpicture
\setcoordinatesystem units <2.5mm,2.5mm>
%\setplotarea x from 0.5 to 1, y from 0 to 3.5
\put {\tiny$#1$} at 0.5 0
\put {\tiny$#2$} at 0.5 1
\put {\tiny$#3$} at 0.5 2
\put {\tiny$#4$} at 2   1.6
\put {$\cdot$} at 1 0
\put {$\cdot$} at 1 1
\put {$\cdot$} at 1 2
\put {$\cdot$} at 2 1
\setlinear
\plot 1 0 1 2 /
\plot 1 0 2 1 /
\endpicture
}

\newcommand{\tY}[4]
{
\beginpicture
\setcoordinatesystem units <2.5mm,2.5mm>
%\setplotarea x from 0.5 to 1, y from 0 to 3.8
\put {\tiny$#1$} at 0.5 0
\put {\tiny$#2$} at 0.5 1
\put {\tiny$#3$} at 0.5 2.6
\put {\tiny$#4$} at 1.5 2.6
\put {$\cdot$} at 1 0
\put {$\cdot$} at 1 1
\put {$\cdot$} at 0.5 2
\put {$\cdot$} at 1.5 2
\setlinear
\plot 1 0 1 1 0.5 2 /
\plot 1 1 1.5 2 /
\endpicture
}

\newcommand{\tI}[4]
{
\beginpicture
\setcoordinatesystem units <2.5mm,2.5mm>
%\setplotarea x from 0.5 to 1, y from 0 to 3.5
\put {\tiny$#1$} at 0.5 0
\put {\tiny$#2$} at 0.5 1
\put {\tiny$#3$} at 0.5 2
\put {\tiny$#4$} at 0.5 3
\put {$\cdot$} at 1 0
\put {$\cdot$} at 1 1
\put {$\cdot$} at 1 2
\put {$\cdot$} at 1 3
\setlinear
\plot 1 0 1 3 /
\endpicture
}

\begin{table}
\caption{Labeled trees for $P_4^2$}
\label{table2P42}
\[
\begin{tabular}{|cccc|r@{$\ =\ $}r@{$\;\cdot\:4!\ $}|r|l|}
\hline
\multicolumn{4}{|c|}{Labeled Trees} &
\multicolumn{2}{c|}{\# (b)} &
\# (d) &
Example \\
\hline
\tW{a}{b}{c}{d} & \tW{b}{a}{c}{d} & \tW{c}{a}{b}{d} & \tW{d}{a}{b}{c} &
24 & 1 & 2160
 & \scriptsize{Id} \\
\hline
\tW{a}{b}{c}{d} & \tW{b}{a}{c}{d} & \tL{c}{d}{b}{a} & \tL{d}{c}{a}{b} &
288 & 12 & 576
 & \scriptsize{(4, 7)} \\
\tW{a}{b}{c}{d} & \tW{b}{a}{c}{d} & \tL{c}{d}{b}{a} & \tL{d}{c}{b}{a} &
288 & 12 & 2160
 & \scriptsize{(7, 8)} \\
\hline
\tW{a}{b}{c}{d} & \tW{b}{a}{c}{d} & \tY{c}{d}{a}{b} & \tY{d}{c}{a}{b} &
144 & 6 & 2160
 & \scriptsize{(3, 4)(7, 8)} \\
\tW{a}{b}{c}{d} & \tW{b}{a}{c}{d} & \tY{d}{c}{a}{b} & \tY{d}{c}{a}{b} &
144 & 6 & 0
 &  \\
\hline
\tW{a}{b}{c}{d} & \tL{b}{c}{d}{a} & \tL{c}{d}{a}{b} & \tL{d}{b}{a}{c} &
576 & 24 & 576
 & \scriptsize{(3, 7, 14, 10, 6, 4)(5, 13, 9)(8, 16, 12)(11, 15)} \\
\tW{a}{b}{c}{d} & \tL{b}{d}{a}{c} & \tL{c}{b}{a}{d} & \tL{d}{c}{a}{b} &
192 & 8 & 2160
 & \scriptsize{(2, 3, 4)} \\
\hline
\tW{a}{b}{c}{d} & \tL{b}{c}{d}{a} & \tL{d}{c}{a}{b} & \tI{c}{b}{d}{a} &
576 & 24 & 2160
 & \scriptsize{(7, 8)(11, 12)(14, 16, 15)} \\
\tW{a}{b}{c}{d} & \tL{c}{b}{a}{d} & \tL{d}{b}{a}{c} & \tI{b}{d}{c}{a} &
576 & 24 & 576
 & \scriptsize{(3, 4, 7, 10, 8, 11, 6)(5, 9)(14, 16, 15)} \\
\hline
\tW{a}{b}{c}{d} & \tL{d}{b}{a}{c} & \tY{c}{b}{a}{d} & \tI{c}{b}{d}{a} &
576 & 24 & 0
 &  \\
\tW{a}{b}{c}{d} & \tL{d}{c}{a}{b} & \tY{c}{b}{a}{d} & \tI{b}{c}{d}{a} &
576 & 24 & 2160
 & \scriptsize{(2, 3)(6, 8)(10, 12)(14, 15, 16)} \\
\hline
\tL{a}{b}{d}{c} & \tL{b}{a}{c}{d} & \tY{c}{d}{a}{b} & \tY{d}{c}{a}{b} &
288 & 12 & 576
 & \scriptsize{(3, 4)(7, 8)(10, 13)} \\
\tL{a}{b}{d}{c} & \tL{b}{a}{c}{d} & \tY{d}{c}{a}{b} & \tY{d}{c}{a}{b} &
288 & 12 & 0
 &  \\
\tL{a}{b}{d}{c} & \tL{b}{a}{d}{c} & \tY{c}{d}{a}{b} & \tY{c}{d}{a}{b} &
144 & 6 & 0
 &  \\
\tL{a}{b}{d}{c} & \tL{b}{a}{d}{c} & \tY{c}{d}{a}{b} & \tY{d}{c}{a}{b} &
288 & 12 & 2160
 & \scriptsize{(3, 4)(7, 8)(13, 14)} \\
\tL{a}{b}{d}{c} & \tL{b}{a}{d}{c} & \tY{d}{c}{a}{b} & \tY{d}{c}{a}{b} &
144 & 6 & 0
 &  \\
\hline
\tY{a}{b}{c}{d} & \tY{a}{b}{c}{d} & \tY{c}{d}{a}{b} & \tY{d}{c}{a}{b} &
144 & 6 & 0
 &  \\
\tY{a}{b}{c}{d} & \tY{a}{b}{c}{d} & \tY{d}{c}{a}{b} & \tY{d}{c}{a}{b} &
72 & 3 & 1536
 & \scriptsize{(2, 9, 5)(4, 11, 7)(6, 10, 13)(8, 12, 15)} \\
\tY{a}{b}{c}{d} & \tY{b}{a}{c}{d} & \tY{c}{d}{a}{b} & \tY{d}{c}{a}{b} &
72 & 3 & 2160
 & \scriptsize{(3, 4)(7, 8)(9, 10)(13, 14)} \\
\hline
\end{tabular}
\]
\end{table}

\section{Additional examples}\label{section3}

We wish to relate the above analysis to the automorphisms constructed in \cite{Sz},
namely Examples 8 and 9 therein.

\begin{example}\label{psi}\rm
Consider a nontrivial partition $W_n^1 = R_1 \cup \ldots \cup R_r$ of $W_n^1$ into a union of $r$ disjoint subsets,
$1 < r \leq n$.
Let $\sigma_i \in P_n^1$, $i=1,\ldots,r$, be permutations of $W_n^1$ such that $\sigma_i \sigma_j^{-1}(R_m) = R_m$
for all $i,j,m \in \{1,\ldots,r\}$. We define $\psi \in P_n^2$ as
$\psi(\alpha,\beta) = (\alpha,\sigma_i(\beta))$
for $\alpha \in R_i$, $\beta \in W_n^1$.
So constructed $\lambda_\psi$ is invertible, with inverse $\lambda_{\overline\psi}$ where $\overline\psi \in P_n^3$
is given by
$\overline\psi(\alpha,\beta,\gamma) = (\alpha,\sigma_i^{-1}(\beta),\sigma_j \sigma_k^{-1}(\gamma))$
for $\alpha \in R_i$, $\beta \in R_k$, $\sigma_i^{-1}(\beta) \in R_j$.
Moreover, it is easy to see that $\lambda_\psi \in {\rm Inn}(\O_n)$ if and only if $\psi = {\rm id}$.
\end{example}

If $n=4$, $r=2$, $R_1 = \{1,2\}$, $R_2 = \{3,4\}$, $\sigma_1 = (23) (= \sigma_1^{-1})$, $\sigma_2 = (1243)$,
$\psi$ is constructed  from these data as above and $w = S_1 S_1^* + S_3 S_2^* + S_2 S_3^* + S_4 S_4^*  \sim \sigma_1$,
then ${\rm Ad}(w) \lambda_\psi$ is the outer automorphism of $\O_4$ constructed and discussed
by Matsumoto and Tomiyama in \cite{MaTo}.
% Also, by a similar argument for $n=2=r$ one recovers the flip-flop automorphism of $\O_2$
For this specific example, it is not difficult to verify that
the corresponding 4-tuple of rooted trees is of type $\alpha\alpha\alpha\alpha$
according to Table \ref{table1P42}.

More generally, for any $\lambda_\psi \in {\rm Aut}(\O_n)$ constructed as in
Example \ref{psi}, the corresponding $n$-tuple of rooted trees can be easily
described as follows. Each tree has $n$ vertices labeled by the elements in
$W_n^1 = \{1,\ldots,n\}$,
and the $i$th tree has root $i$ and all the other
vertices are connected to the root. This readily follows from the fact that the
defining relation (see \cite[Section 4.1]{CS})
$$(i,\alpha) = \psi(\beta,m), \quad \alpha,\beta \in W_n^1$$
for some $m \in \{i,\ldots,n\}$ forces $\beta = i$ and then it can
be solved for all $\alpha$'s.

\begin{example}\label{z2z3}\rm
Let $n \geq 3$, $\phi = (123) \in P_n^1$ and let $\psi \in P_n^2$ be constructed as in Example \ref{psi} with $r=2$
from the data: $R_1 = \{1,2\}$, $R_2 = \{3,\ldots,n\}$, $\sigma_1 = {\rm id}$, $\sigma_2 = (12)$.
Then one checks that $\lambda_\phi$ and $\lambda_\psi$ are outer automorphisms of $\O_n$ of order 3 and 2, respectively.
% $\lambda_\phi$ is Bogolubov
and the group generated by $\lambda_\phi$ and $\lambda_\psi$ in ${\rm Out}(\O_n)$ is ${\mathbb Z}_3 * {\mathbb Z}_2$ \cite{Sz}.
\end{example}
Since $\lambda_\phi$ is Bogolubov automorphism, the trees associated to $\phi$
(thought of as an element in $P_n^k$ for $k>1$) are computed in \cite[Example 4.4]{CS}.
Also, as discussed above,
the $n$ trees corresponding to $\psi \in P_n^2$ are also all identical,
with the root receiving $n-1$ edges from the other vertices.

%\clearpage
\section*{Appendix}
\begin{table}[hb]
\caption{Type $A,A,A$, type $A,H,H$ and type $A,F,G$.}
\begin{center}
\begin{tiny}
\begin{tabular}{cc}
\multirow{1}*{
\begin{tabular}{|r|r|r|r@{$\;\cdot\:9!\ $}|r|}
\hline
A & A & A & \multicolumn{1}{r|}{\# (b)} & count \\
\hline\hline
1 & 1 & 1 & 2 & 1 \\
1 & 1 & 1 & 6 & 487 \\
1 & 1 & 2 & 3 & 16 \\
1 & 1 & 2 & 6 & 1032 \\
1 & 2 & 2 & 6 & 568 \\
2 & 2 & 2 & 1 & 1 \\
2 & 2 & 2 & 3 & 9 \\
2 & 2 & 2 & 6 & 54 \\
\hline
\end{tabular}
}
& \multirow{3}*{
\begin{tabular}{|r|r|r|r@{$\;\cdot\:9!\ $}|r|}
\hline
A & F & G & \multicolumn{1}{r|}{\# (b)} & count \\
\hline\hline
1 & 6 & 7 & 6 & 6 \\
1 & 6 & 31 & 6 & 18 \\
1 & 7 & 7 & 6 & 2 \\
1 & 9 & 7 & 6 & 2 \\
2 & 1 & 1 & 6 & 8 \\
2 & 1 & 2 & 6 & 8 \\
2 & 1 & 3 & 6 & 4 \\
2 & 2 & 1 & 6 & 8 \\
2 & 2 & 2 & 6 & 4 \\
2 & 3 & 11 & 6 & 4 \\
2 & 3 & 26 & 6 & 8 \\
2 & 3 & 31 & 6 & 4 \\
2 & 3 & 32 & 6 & 8 \\
2 & 4 & 11 & 6 & 16 \\
2 & 4 & 26 & 6 & 16 \\
2 & 4 & 28 & 6 & 16 \\
2 & 4 & 31 & 6 & 8 \\
2 & 4 & 32 & 6 & 16 \\
2 & 6 & 7 & 6 & 20 \\
2 & 6 & 11 & 6 & 8 \\
2 & 6 & 21 & 6 & 16 \\
2 & 6 & 28 & 6 & 32 \\
2 & 6 & 31 & 6 & 56 \\
2 & 6 & 32 & 6 & 16 \\
2 & 7 & 7 & 6 & 20 \\
2 & 7 & 11 & 6 & 8 \\
2 & 8 & 11 & 6 & 8 \\
2 & 9 & 7 & 6 & 20 \\
2 & 9 & 11 & 6 & 8 \\
2 & 9 & 26 & 6 & 16 \\
2 & 9 & 31 & 6 & 8 \\
\hline
\end{tabular}
}
\\
\\
\\
\\
\\
\\
\\
\\
\\
\\
\\
\\
\\
\\
\\
\phantom{jsk} & \\
\multirow{1}*{
\begin{tabular}{|r|r|r|r@{$\;\cdot\:9!\ $}|r|}
\hline
A & H & H & \multicolumn{1}{r|}{\# (b)} & count \\
\hline\hline
2 & 1 & 11 & 6 & 2 \\
2 & 2 & 9 & 6 & 2 \\
2 & 2 & 11 & 6 & 2 \\
2 & 4 & 5 & 6 & 2 \\
2 & 4 & 9 & 6 & 4 \\
2 & 4 & 11 & 6 & 2 \\
2 & 18 & 23 & 6 & 8 \\
2 & 23 & 23 & 3 & 2 \\
2 & 23 & 23 & 6 & 2 \\
\hline
\end{tabular}
}&\\
\end{tabular}
\end{tiny}
\end{center}
\end{table}

\clearpage
\begin{center}
Table 9: Types $A,C,D$.\par
\bigskip
\begin{tiny}
\begin{tabular}{|r|r|r|r@{$\;\cdot\:9!\ $}|r|}
\hline
A & C & D & \multicolumn{1}{r|}{\# (b)} & count \\
\hline\hline
1 & 1 & 8 & 6 & 12 \\
1 & 1 & 12 & 6 & 12 \\
1 & 5 & 8 & 6 & 6 \\
1 & 7 & 8 & 6 & 6 \\
2 & 1 & 7 & 6 & 8 \\
2 & 1 & 8 & 6 & 36 \\
2 & 1 & 12 & 6 & 36 \\
2 & 2 & 4 & 6 & 8 \\
2 & 2 & 6 & 6 & 8 \\
2 & 2 & 7 & 6 & 10 \\
2 & 2 & 8 & 6 & 12 \\
2 & 2 & 9 & 6 & 10 \\
2 & 2 & 10 & 6 & 4 \\
2 & 2 & 11 & 6 & 36 \\
2 & 2 & 13 & 6 & 8 \\
2 & 3 & 1 & 6 & 12 \\
2 & 3 & 2 & 6 & 12 \\
2 & 3 & 3 & 6 & 8 \\
2 & 3 & 4 & 6 & 6 \\
2 & 3 & 5 & 6 & 18 \\
2 & 3 & 6 & 6 & 10 \\
2 & 3 & 7 & 6 & 8 \\
2 & 3 & 8 & 6 & 6 \\
2 & 3 & 9 & 6 & 6 \\
2 & 3 & 10 & 6 & 5 \\
2 & 3 & 11 & 6 & 30 \\
2 & 3 & 12 & 6 & 3 \\
2 & 3 & 13 & 6 & 9 \\
2 & 3 & 14 & 6 & 6 \\
2 & 4 & 4 & 6 & 2 \\
2 & 4 & 7 & 6 & 6 \\
2 & 4 & 11 & 6 & 18 \\
2 & 5 & 1 & 6 & 18 \\
2 & 5 & 2 & 6 & 6 \\
2 & 5 & 4 & 6 & 6 \\
2 & 5 & 7 & 6 & 14 \\
2 & 5 & 8 & 6 & 16 \\
2 & 5 & 9 & 6 & 3 \\
2 & 5 & 11 & 6 & 40 \\
2 & 5 & 12 & 6 & 6 \\
2 & 5 & 13 & 6 & 4 \\
2 & 6 & 1 & 6 & 18 \\
2 & 6 & 2 & 6 & 18 \\
2 & 6 & 3 & 6 & 10 \\
2 & 6 & 4 & 6 & 12 \\
2 & 6 & 5 & 6 & 8 \\
2 & 6 & 6 & 6 & 4 \\
2 & 6 & 7 & 6 & 12 \\
2 & 6 & 8 & 6 & 12 \\
2 & 6 & 9 & 6 & 12 \\
2 & 6 & 10 & 6 & 10 \\
2 & 6 & 11 & 6 & 56 \\
2 & 6 & 13 & 6 & 8 \\
2 & 6 & 14 & 6 & 3 \\
2 & 7 & 1 & 6 & 6 \\
2 & 7 & 2 & 6 & 3 \\
2 & 7 & 4 & 6 & 6 \\
2 & 7 & 7 & 6 & 14 \\
2 & 7 & 8 & 6 & 16 \\
2 & 7 & 9 & 6 & 3 \\
2 & 7 & 11 & 6 & 38 \\
2 & 7 & 13 & 6 & 2 \\
2 & 8 & 1 & 6 & 24 \\
2 & 8 & 2 & 6 & 15 \\
2 & 8 & 3 & 6 & 5 \\
2 & 8 & 4 & 6 & 18 \\
2 & 8 & 5 & 6 & 9 \\
2 & 8 & 6 & 6 & 4 \\
2 & 8 & 7 & 6 & 26 \\
2 & 8 & 8 & 6 & 6 \\
2 & 8 & 9 & 6 & 3 \\
2 & 8 & 10 & 6 & 2 \\
2 & 8 & 11 & 6 & 66 \\
2 & 8 & 13 & 6 & 8 \\
2 & 8 & 14 & 6 & 3 \\
\hline
\end{tabular}
\end{tiny}
\end{center}
%\end{table}

\begin{table}
\caption{Type $A,B,B$.}
\medskip
\begin{tiny}
\begin{center}
\begin{tabular}{ccc}
\begin{tabular}{|r|r|r|r@{$\;\cdot\:9!\ $}|r|}
\hline
A & B & B & \multicolumn{1}{r|}{\# (b)} & count \\
\hline\hline
1 & 2 & 2 & 6 & 2 \\
1 & 2 & 3 & 6 & 16 \\
1 & 2 & 4 & 6 & 12 \\
1 & 2 & 5 & 6 & 4 \\
1 & 2 & 7 & 6 & 8 \\
1 & 2 & 8 & 6 & 30 \\
1 & 2 & 9 & 6 & 12 \\
1 & 2 & 11 & 6 & 28 \\
1 & 2 & 12 & 6 & 28 \\
1 & 2 & 13 & 6 & 12 \\
1 & 2 & 14 & 6 & 16 \\
1 & 2 & 15 & 6 & 4 \\
1 & 2 & 16 & 6 & 20 \\
1 & 2 & 18 & 6 & 14 \\
1 & 3 & 7 & 6 & 12 \\
1 & 3 & 15 & 6 & 6 \\
1 & 4 & 4 & 6 & 2 \\
1 & 4 & 5 & 6 & 4 \\
1 & 4 & 7 & 6 & 8 \\
1 & 4 & 9 & 6 & 4 \\
1 & 4 & 13 & 6 & 4 \\
1 & 4 & 15 & 6 & 4 \\
1 & 4 & 18 & 6 & 4 \\
1 & 5 & 7 & 6 & 8 \\
1 & 5 & 9 & 6 & 4 \\
1 & 5 & 13 & 6 & 4 \\
1 & 5 & 15 & 6 & 4 \\
1 & 5 & 18 & 6 & 4 \\
1 & 7 & 7 & 6 & 6 \\
1 & 7 & 8 & 6 & 18 \\
1 & 7 & 9 & 6 & 4 \\
1 & 7 & 11 & 6 & 12 \\
1 & 7 & 12 & 6 & 12 \\
1 & 7 & 13 & 6 & 8 \\
1 & 7 & 15 & 6 & 4 \\
1 & 7 & 16 & 6 & 28 \\
1 & 7 & 18 & 6 & 18 \\
1 & 8 & 15 & 6 & 8 \\
1 & 9 & 13 & 6 & 4 \\
1 & 9 & 18 & 6 & 4 \\
1 & 11 & 15 & 6 & 8 \\
1 & 12 & 15 & 6 & 12 \\
1 & 13 & 18 & 6 & 4 \\
1 & 14 & 15 & 6 & 6 \\
1 & 15 & 16 & 6 & 8 \\
1 & 15 & 18 & 6 & 4 \\
1 & 18 & 18 & 6 & 2 \\
2 & 1 & 12 & 6 & 12 \\
2 & 1 & 14 & 6 & 6 \\
2 & 1 & 16 & 6 & 6 \\
2 & 2 & 2 & 3 & 1 \\
2 & 2 & 2 & 6 & 2 \\
2 & 2 & 3 & 6 & 24 \\
2 & 2 & 4 & 6 & 18 \\
2 & 2 & 5 & 6 & 6 \\
2 & 2 & 6 & 6 & 4 \\
2 & 2 & 7 & 6 & 12 \\
2 & 2 & 8 & 6 & 48 \\
2 & 2 & 9 & 6 & 18 \\
2 & 2 & 10 & 6 & 12 \\
2 & 2 & 11 & 6 & 48 \\
2 & 2 & 12 & 6 & 52 \\
2 & 2 & 13 & 6 & 18 \\
2 & 2 & 14 & 6 & 26 \\
2 & 2 & 15 & 6 & 6 \\
2 & 2 & 16 & 6 & 26 \\
2 & 2 & 17 & 6 & 4 \\
2 & 2 & 18 & 6 & 18 \\
\hline
\end{tabular}
&
\begin{tabular}{|r|r|r|r@{$\;\cdot\:9!\ $}|r|}
\hline
A & B & B & \multicolumn{1}{r|}{\# (b)} & count \\
\hline\hline
2 & 3 & 3 & 3 & 1 \\
2 & 3 & 3 & 6 & 10 \\
2 & 3 & 4 & 6 & 8 \\
2 & 3 & 6 & 6 & 2 \\
2 & 3 & 7 & 6 & 24 \\
2 & 3 & 8 & 6 & 44 \\
2 & 3 & 9 & 6 & 8 \\
2 & 3 & 10 & 6 & 20 \\
2 & 3 & 11 & 6 & 24 \\
2 & 3 & 12 & 6 & 32 \\
2 & 3 & 13 & 6 & 8 \\
2 & 3 & 14 & 6 & 16 \\
2 & 3 & 15 & 6 & 12 \\
2 & 3 & 16 & 6 & 44 \\
2 & 3 & 17 & 6 & 4 \\
2 & 3 & 18 & 6 & 24 \\
2 & 4 & 4 & 3 & 1 \\
2 & 4 & 4 & 6 & 2 \\
2 & 4 & 5 & 6 & 6 \\
2 & 4 & 7 & 6 & 12 \\
2 & 4 & 8 & 6 & 16 \\
2 & 4 & 9 & 6 & 6 \\
2 & 4 & 10 & 6 & 16 \\
2 & 4 & 12 & 6 & 20 \\
2 & 4 & 13 & 6 & 6 \\
2 & 4 & 14 & 6 & 8 \\
2 & 4 & 15 & 6 & 6 \\
2 & 4 & 16 & 6 & 18 \\
2 & 4 & 17 & 6 & 4 \\
2 & 4 & 18 & 6 & 14 \\
2 & 5 & 7 & 6 & 12 \\
2 & 5 & 9 & 6 & 6 \\
2 & 5 & 12 & 6 & 8 \\
2 & 5 & 13 & 6 & 6 \\
2 & 5 & 15 & 6 & 6 \\
2 & 5 & 18 & 6 & 6 \\
2 & 6 & 6 & 3 & 1 \\
2 & 6 & 6 & 6 & 4 \\
2 & 6 & 7 & 6 & 8 \\
2 & 6 & 8 & 6 & 4 \\
2 & 6 & 10 & 6 & 6 \\
2 & 6 & 12 & 6 & 16 \\
2 & 6 & 14 & 6 & 16 \\
2 & 6 & 15 & 6 & 6 \\
2 & 6 & 16 & 6 & 10 \\
2 & 6 & 17 & 6 & 4 \\
2 & 7 & 7 & 3 & 1 \\
2 & 7 & 7 & 6 & 8 \\
2 & 7 & 8 & 6 & 48 \\
2 & 7 & 9 & 6 & 6 \\
2 & 7 & 10 & 6 & 16 \\
2 & 7 & 11 & 6 & 24 \\
2 & 7 & 12 & 6 & 32 \\
2 & 7 & 13 & 6 & 12 \\
2 & 7 & 14 & 6 & 8 \\
2 & 7 & 15 & 6 & 6 \\
2 & 7 & 16 & 6 & 40 \\
2 & 7 & 18 & 6 & 24 \\
\hline
\end{tabular}
&
\begin{tabular}{|r|r|r|r@{$\;\cdot\:9!\ $}|r|}
\hline
A & B & B & \multicolumn{1}{r|}{\# (b)} & count \\
\hline\hline
2 & 8 & 8 & 3 & 2 \\
2 & 8 & 8 & 6 & 42 \\
2 & 8 & 9 & 6 & 16 \\
2 & 8 & 10 & 6 & 36 \\
2 & 8 & 11 & 6 & 48 \\
2 & 8 & 12 & 6 & 56 \\
2 & 8 & 13 & 6 & 16 \\
2 & 8 & 14 & 6 & 26 \\
2 & 8 & 15 & 6 & 24 \\
2 & 8 & 16 & 6 & 88 \\
2 & 8 & 17 & 6 & 4 \\
2 & 8 & 18 & 6 & 48 \\
2 & 9 & 10 & 6 & 16 \\
2 & 9 & 12 & 6 & 20 \\
2 & 9 & 13 & 6 & 6 \\
2 & 9 & 14 & 6 & 8 \\
2 & 9 & 16 & 6 & 12 \\
2 & 9 & 17 & 6 & 4 \\
2 & 9 & 18 & 6 & 8 \\
2 & 10 & 11 & 6 & 32 \\
2 & 10 & 12 & 6 & 44 \\
2 & 10 & 13 & 6 & 16 \\
2 & 10 & 14 & 6 & 26 \\
2 & 10 & 15 & 6 & 12 \\
2 & 10 & 16 & 6 & 34 \\
2 & 10 & 18 & 6 & 20 \\
2 & 11 & 12 & 6 & 8 \\
2 & 11 & 14 & 6 & 12 \\
2 & 11 & 15 & 6 & 12 \\
2 & 11 & 16 & 6 & 54 \\
2 & 11 & 17 & 6 & 12 \\
2 & 11 & 18 & 6 & 24 \\
2 & 12 & 12 & 6 & 8 \\
2 & 12 & 13 & 6 & 24 \\
2 & 12 & 14 & 6 & 28 \\
2 & 12 & 15 & 6 & 28 \\
2 & 12 & 16 & 6 & 64 \\
2 & 12 & 17 & 6 & 20 \\
2 & 12 & 18 & 6 & 44 \\
2 & 13 & 14 & 6 & 4 \\
2 & 13 & 16 & 6 & 12 \\
2 & 13 & 17 & 6 & 4 \\
2 & 13 & 18 & 6 & 8 \\
2 & 14 & 14 & 3 & 2 \\
2 & 14 & 14 & 6 & 8 \\
2 & 14 & 15 & 6 & 14 \\
2 & 14 & 16 & 6 & 24 \\
2 & 14 & 17 & 6 & 10 \\
2 & 14 & 18 & 6 & 12 \\
2 & 15 & 16 & 6 & 14 \\
2 & 15 & 17 & 6 & 4 \\
2 & 15 & 18 & 6 & 6 \\
2 & 16 & 16 & 3 & 3 \\
2 & 16 & 16 & 6 & 26 \\
2 & 16 & 17 & 6 & 10 \\
2 & 16 & 18 & 6 & 30 \\
2 & 17 & 18 & 6 & 4 \\
2 & 18 & 18 & 3 & 2 \\
2 & 18 & 18 & 6 & 8 \\
\hline
\end{tabular}
\end{tabular}
\end{center}
\end{tiny}
\end{table}

\begin{table}
\caption{Type $A,E,E$.}
\medskip
\begin{tiny}
\begin{tabular}{ccc}
\begin{tabular}{|r|r|r|r@{$\;\cdot\:9!\ $}|r|}
\hline
A & E & E & \multicolumn{1}{r|}{\# (b)} & count \\
\hline\hline
1 & 1 & 3 & 6 & 2 \\
1 & 1 & 5 & 6 & 2 \\
1 & 13 & 13 & 6 & 4 \\
2 & 1 & 3 & 6 & 6 \\
2 & 1 & 5 & 6 & 6 \\
2 & 2 & 3 & 6 & 4 \\
2 & 2 & 4 & 6 & 2 \\
2 & 3 & 23 & 6 & 4 \\
2 & 3 & 38 & 6 & 8 \\
2 & 3 & 43 & 6 & 4 \\
2 & 6 & 22 & 6 & 2 \\
2 & 6 & 34 & 6 & 2 \\
2 & 6 & 41 & 6 & 2 \\
2 & 6 & 42 & 6 & 2 \\
2 & 8 & 8 & 3 & 1 \\
2 & 8 & 20 & 6 & 2 \\
2 & 8 & 22 & 6 & 4 \\
2 & 8 & 23 & 6 & 2 \\
2 & 8 & 27 & 6 & 2 \\
2 & 8 & 32 & 6 & 2 \\
2 & 8 & 34 & 6 & 4 \\
2 & 8 & 35 & 6 & 2 \\
2 & 8 & 36 & 6 & 2 \\
2 & 8 & 38 & 6 & 4 \\
2 & 8 & 40 & 6 & 2 \\
2 & 8 & 41 & 6 & 4 \\
2 & 8 & 42 & 6 & 4 \\
2 & 8 & 43 & 6 & 2 \\
2 & 8 & 45 & 6 & 2 \\
2 & 9 & 20 & 6 & 8 \\
2 & 9 & 23 & 6 & 16 \\
2 & 9 & 32 & 6 & 8 \\
2 & 9 & 35 & 6 & 8 \\
2 & 9 & 36 & 6 & 8 \\
2 & 9 & 38 & 6 & 16 \\
2 & 9 & 43 & 6 & 8 \\
2 & 9 & 44 & 6 & 8 \\
2 & 9 & 45 & 6 & 8 \\
2 & 10 & 20 & 6 & 8 \\
2 & 10 & 23 & 6 & 4 \\
2 & 10 & 32 & 6 & 8 \\
2 & 10 & 35 & 6 & 8 \\
2 & 10 & 36 & 6 & 8 \\
2 & 10 & 38 & 6 & 8 \\
2 & 10 & 43 & 6 & 4 \\
2 & 10 & 44 & 6 & 4 \\
2 & 10 & 45 & 6 & 8 \\
2 & 11 & 20 & 6 & 12 \\
2 & 11 & 32 & 6 & 8 \\
2 & 11 & 35 & 6 & 4 \\
2 & 11 & 36 & 6 & 8 \\
2 & 11 & 45 & 6 & 8 \\
2 & 12 & 12 & 3 & 2 \\
2 & 12 & 12 & 6 & 10 \\
2 & 12 & 14 & 6 & 16 \\
2 & 12 & 18 & 6 & 8 \\
2 & 12 & 21 & 6 & 32 \\
2 & 12 & 23 & 6 & 8 \\
2 & 12 & 31 & 6 & 16 \\
2 & 12 & 33 & 6 & 32 \\
2 & 12 & 37 & 6 & 32 \\
2 & 12 & 38 & 6 & 16 \\
2 & 12 & 43 & 6 & 16 \\
2 & 12 & 46 & 6 & 32 \\
\hline
\end{tabular}
&
\begin{tabular}{|r|r|r|r@{$\;\cdot\:9!\ $}|r|}
\hline
A & E & E & \multicolumn{1}{r|}{\# (b)} & count \\
\hline\hline
2 & 13 & 13 & 3 & 2 \\
2 & 13 & 13 & 6 & 10 \\
2 & 13 & 23 & 6 & 8 \\
2 & 15 & 20 & 6 & 2 \\
2 & 15 & 22 & 6 & 4 \\
2 & 15 & 27 & 6 & 2 \\
2 & 15 & 32 & 6 & 2 \\
2 & 15 & 34 & 6 & 2 \\
2 & 15 & 36 & 6 & 2 \\
2 & 15 & 41 & 6 & 2 \\
2 & 15 & 42 & 6 & 2 \\
2 & 15 & 45 & 6 & 2 \\
2 & 17 & 20 & 6 & 12 \\
2 & 17 & 32 & 6 & 4 \\
2 & 17 & 36 & 6 & 4 \\
2 & 17 & 45 & 6 & 4 \\
2 & 18 & 43 & 6 & 8 \\
2 & 18 & 46 & 6 & 16 \\
2 & 19 & 20 & 6 & 4 \\
2 & 20 & 20 & 3 & 1 \\
2 & 20 & 20 & 6 & 2 \\
2 & 20 & 22 & 6 & 14 \\
2 & 20 & 25 & 6 & 2 \\
2 & 20 & 27 & 6 & 4 \\
2 & 20 & 28 & 6 & 4 \\
2 & 20 & 30 & 6 & 24 \\
2 & 20 & 32 & 6 & 6 \\
2 & 20 & 34 & 6 & 12 \\
2 & 20 & 36 & 6 & 6 \\
2 & 20 & 39 & 6 & 2 \\
2 & 20 & 41 & 6 & 12 \\
2 & 20 & 42 & 6 & 12 \\
2 & 20 & 45 & 6 & 6 \\
2 & 22 & 22 & 3 & 1 \\
2 & 22 & 22 & 6 & 2 \\
2 & 22 & 23 & 6 & 2 \\
2 & 22 & 24 & 6 & 4 \\
2 & 22 & 25 & 6 & 6 \\
2 & 22 & 27 & 6 & 4 \\
2 & 22 & 28 & 6 & 8 \\
2 & 22 & 32 & 6 & 6 \\
2 & 22 & 34 & 6 & 4 \\
2 & 22 & 35 & 6 & 4 \\
2 & 22 & 36 & 6 & 6 \\
2 & 22 & 38 & 6 & 4 \\
2 & 22 & 39 & 6 & 4 \\
2 & 22 & 40 & 6 & 2 \\
2 & 22 & 41 & 6 & 4 \\
2 & 22 & 42 & 6 & 4 \\
2 & 22 & 43 & 6 & 2 \\
2 & 22 & 44 & 6 & 2 \\
2 & 22 & 45 & 6 & 6 \\
2 & 23 & 38 & 6 & 2 \\
2 & 23 & 42 & 6 & 2 \\
2 & 23 & 43 & 6 & 8 \\
2 & 23 & 45 & 6 & 2 \\
2 & 23 & 46 & 6 & 16 \\
2 & 24 & 34 & 6 & 2 \\
2 & 24 & 41 & 6 & 2 \\
2 & 24 & 42 & 6 & 4 \\
\hline
\end{tabular}
&
\begin{tabular}{|r|r|r|r@{$\;\cdot\:9!\ $}|r|}
\hline
A & E & E & \multicolumn{1}{r|}{\# (b)} & count \\
\hline\hline
2 & 25 & 32 & 6 & 2 \\
2 & 25 & 34 & 6 & 4 \\
2 & 25 & 36 & 6 & 2 \\
2 & 25 & 41 & 6 & 4 \\
2 & 25 & 42 & 6 & 4 \\
2 & 25 & 45 & 6 & 2 \\
2 & 27 & 28 & 6 & 4 \\
2 & 27 & 32 & 6 & 4 \\
2 & 27 & 34 & 6 & 2 \\
2 & 27 & 36 & 6 & 4 \\
2 & 27 & 38 & 6 & 2 \\
2 & 27 & 39 & 6 & 2 \\
2 & 27 & 41 & 6 & 2 \\
2 & 27 & 42 & 6 & 4 \\
2 & 27 & 45 & 6 & 6 \\
2 & 28 & 32 & 6 & 2 \\
2 & 28 & 34 & 6 & 2 \\
2 & 28 & 36 & 6 & 4 \\
2 & 28 & 38 & 6 & 2 \\
2 & 28 & 41 & 6 & 4 \\
2 & 28 & 42 & 6 & 2 \\
2 & 28 & 45 & 6 & 2 \\
2 & 30 & 32 & 6 & 8 \\
2 & 30 & 36 & 6 & 8 \\
2 & 30 & 45 & 6 & 8 \\
2 & 31 & 43 & 6 & 16 \\
2 & 31 & 46 & 6 & 32 \\
2 & 32 & 36 & 6 & 6 \\
2 & 32 & 38 & 6 & 4 \\
2 & 32 & 41 & 6 & 12 \\
2 & 32 & 45 & 6 & 2 \\
2 & 34 & 35 & 6 & 2 \\
2 & 34 & 36 & 6 & 6 \\
2 & 34 & 38 & 6 & 4 \\
2 & 34 & 39 & 6 & 2 \\
2 & 34 & 40 & 6 & 2 \\
2 & 34 & 41 & 6 & 2 \\
2 & 34 & 42 & 6 & 2 \\
2 & 34 & 45 & 6 & 2 \\
2 & 35 & 42 & 6 & 2 \\
2 & 36 & 36 & 3 & 1 \\
2 & 36 & 36 & 6 & 2 \\
2 & 36 & 39 & 6 & 2 \\
2 & 36 & 41 & 6 & 4 \\
2 & 36 & 42 & 6 & 2 \\
2 & 36 & 45 & 6 & 4 \\
2 & 38 & 38 & 3 & 1 \\
2 & 38 & 42 & 6 & 2 \\
2 & 38 & 45 & 6 & 4 \\
2 & 38 & 46 & 6 & 16 \\
2 & 39 & 41 & 6 & 2 \\
2 & 40 & 41 & 6 & 2 \\
2 & 40 & 42 & 6 & 2 \\
2 & 41 & 41 & 3 & 1 \\
2 & 41 & 42 & 6 & 2 \\
2 & 41 & 45 & 6 & 4 \\
2 & 42 & 45 & 6 & 2 \\
2 & 43 & 43 & 3 & 2 \\
2 & 43 & 43 & 6 & 2 \\
2 & 43 & 45 & 6 & 2 \\
2 & 43 & 46 & 6 & 16 \\
2 & 45 & 45 & 6 & 2 \\
\hline
\end{tabular}
\end{tabular}
\end{tiny}
\end{table}

\bigskip

\noindent
Roberto Conti\\
Mathematics, School of Mathematical and Physical Sciences \\
University of Newcastle, Callaghan, NSW 2308, Australia\\
E-mail: Roberto.Conti@newcastle.edu.au \\

\medskip \noindent
Jason Kimberley\\
Mathematics, School of Mathematical and Physical Sciences \\
University of Newcastle, Callaghan, NSW 2308, Australia\\
E-mail: Jason.Kimberley@newcastle.edu.au\\

\medskip \noindent
Wojciech Szyma{\'n}ski\\
Mathematics, School of Mathematical and Physical Sciences \\
University of Newcastle, Callaghan, NSW 2308, Australia \\
% E-mail: Wojciech.Szymanski@newcastle.edu.au
%\medskip
and\\
%\medskip
%\noindent
Department of Mathematics and Computer Science \\
The University of Southern Denmark \\
Campusvej 55, DK-5230 Odense M, Denmark \\
E-mail: szymanski@imada.sdu.dk

\end{document}